\documentclass[12pt,twoside]{article}

\usepackage{amsmath,amssymb}
\usepackage{amsthm}
\usepackage{stmaryrd}
\usepackage{amscd}
\usepackage{fancyhdr}
\usepackage{times}
\usepackage{sectsty}
\usepackage[dvips,pdftex]{geometry}

\sectionfont{\normalsize\bfseries\scshape\centering\MakeUppercase }

\makeatletter
  \def\@seccntformat#1{\csname the#1\endcsname.\quad}
\makeatother

\newcommand{\rest}{\mbox{\parbox[t]{0.1cm}{$|$\\[-10pt] $|$}}}

\newtheoremstyle{fact}
     {\topsep}
     {\topsep}
     {\slshape}
     {}
     {\bfseries}
     {}
     { }
     {\thmname{#1}\thmnumber{ #2.}\thmnote{ \rm (#3)}} 

\newtheorem{theorem}{Theorem}[section]

\newtheorem{proposition}[theorem]{Proposition}
\newtheorem{corollary}[theorem]{Corollary}
\newtheorem{problem}{Problem}

\theoremstyle{definition}
\newtheorem*{remark}{Remark}

\theoremstyle{fact}
\newtheorem{fact}{Fact}

\newtheorem{observation}[theorem]{Observation}

\def\proofont{\fontseries{bx}\fontshape{sc}\selectfont}
\def\proofname{Proof. }

\makeatletter
\renewenvironment{proof}[1][\proofname]{\par
  \normalfont
  \topsep6\p@\@plus6\p@ \trivlist
  \item[\hskip\labelsep\noindent\proofont #1]\ignorespaces
}{%
  \qed\endtrivlist
}
\makeatother

\author{G\'abor Luk\'acs
\thanks{I gratefully acknowledge the generous financial support received 
from the Alexander von Humboldt Foundation that enabled me to do this 
research.}}

\title{Hereditarily $h$-complete groups
}

\begin{document}

\makeatletter
\let\mytitle\@title
\chead{\small\itshape G. Luk\'acs / \mytitle }
\fancyhead[RO,LE]{\small \thepage}
\makeatother

\maketitle 

\def\thanks#1{}

\thispagestyle{empty}

\begin{abstract}

A topological group $G$ is {\em $h$-complete} if every continuous
homomorphic image of $G$ is (Ra\u{\i}kov-)complete; we say that $G$ is
{\em hereditarily $h$-complete} if every closed subgroup of $G$ is
$h$-complete. In this paper, we establish open-map properties of
hereditarily $h$-complete groups with respect to large classes of groups,
and prove a theorem on the (total) minimality of subdirectly represented
groups. Numerous applications are presented, among them: 1. Every 
hereditarily $h$-complete group with quasi-invariant basis is the
projective limit of its metrizable quotients; 2. If every countable 
discrete hereditarily $h$-complete group is finite, then every 
locally compact hereditarily $h$-complete group that has small invariant 
neighborhoods is compact. In the sequel, several open problems are 
formulated.
\end{abstract}

\section*{Introduction}

By the well known Kuratowski-Mr\'{o}wka Theorem, a (Hausdorff) topological 
space $X$ is compact if and only if for any (Hausdorff) topological space 
$Y$ the projection $p_Y: X \times Y \rightarrow Y$ is closed. Inspired by 
this theorem, one says that a topological group $G$ is {\em categorically 
compact} (or briefly,  {\em c-compact}) if for any topological group $H$ 
the image of every  closed subgroup of $G \times H$ under the projection 
$\pi_H: G \times H \rightarrow H$ is closed in $H$. (All topological 
groups in this paper are assumed to be Hausdorff.) The problem of whether 
every $c$-compact topological group is compact has been an open question 
for more than ten years.

The most extensive study of $c$-compact topological groups was done by
Dikranjan and Uspenskij in \cite{DikUsp}, which has also been a source of
inspiration for part of the author's PhD dissertation \cite{GLPHD}.

Categorical compactness (in the general setting of structured
sets) was introduced in 1974 by Manes \cite{Manes} and studied by
Dikranjan and Giuli \cite{DikGiu}, and by Clementino, Giuli and Tholen
\cite{CGT} who consider in particular categorically compact groups.
In \cite{CleTho} Clementino and Tholen proved (independently of
Dikranjan and Uspenskij \cite{DikUsp}, but in greater generality) the
Tychonoff Theorem for $c$-compact groups. 

\pagebreak[3]

The notion of hereditary $h$-completeness we introduce in
Section~\ref{sect:pre} is motivated by the observation that all known
results relating $c$-compactness to any compactness-like property remain
valid if $c$-compactness is replaced by hereditary $h$-completeness.  Our
main result is Theorem~\ref{thm:openmap}, stating that hereditarily
$h$-complete groups satisfy open-map properties with respect to large
classes of groups (groups of countable tightness and $k$-groups); its
immediate corollary generalizes \cite[3.2]{DikUsp}.\linebreak[3] As an
application, we examine groups admitting a quasi-invariant basis (see
Section~\ref{sect:quasi} for definition), and prove a structure theorem
for hereditarily $h$-complete groups with this property
(Theorem~\ref{thm:QIB:struct}). In the course of this investigation, an
important result on (total) minimality of groups subdirectly represented
by (totally) minimal groups is established (Theorem~\ref{thm:subdir:tm}),
which turns out to have several consequences. We conclude with a reduction
theorem (Theorem~\ref{thm:red-SIN:cc}), relating the compactness of
locally compact $c$-compact groups with small invariant neighborhoods to
finiteness of countable discrete $c$-compact groups. In the sequel,
several open problems are formulated.

\section{Preliminaries}

\label{sect:pre}

A topological group is {\em minimal} if it does not admit a coarser  
(Hausdorff) group topology; a group $G$ is {\em totally minimal} if every 
continuous surjective homomorphism $\varphi: G \rightarrow H$ is open, or 
equivalently, if every quotient of $G$ is minimal. It is an
important result  that {\slshape every closed separable 
subgroup of a $c$-compact group is totally minimal} (\cite[3.6]{DikUsp}). 
This is obtained as a consequence of \cite[3.4]{DikUsp}, which requires 
only the closed normal subgroups of the separable subgroup to be $h$-complete, 
and thus the condition of $c$-compactness can be slightly weakened, as we 
explain below.

A topological group $G$ is said to be {\em $h$-complete} if for any 
continuous homomorphism\linebreak $\varphi: G \rightarrow H$ into a 
topological group, $\varphi(G)$ is closed in $H$ (\cite{DikTon}). 
Equivalently, $G$ is $h$-complete if every continuous homomorphic image 
of $G$ is Ra\u{\i}kov-complete (i.e., complete in the uniformity given by 
the join of the left and the right uniformity of the group).
Every $c$-compact group is $h$-complete; moreover, since $c$-compactness
is a closed-hereditary property, every closed subgroup of a $c$-compact
group is $h$-complete. Motivated by this, the author introduces the 
concept of hereditary $h$-completeness in the obvious way: $G$ is {\em  
hereditarily $h$-complete} if every closed subgroup of $G$ is $h$-complete.
The following immediate consequence of \cite[3.4]{DikUsp} will play an 
important role in establishing the results of this paper.

\begin{observation} \label{obs:sepsg:tm}
Every closed separable subgroup of a hereditarily $h$-complete group is 
totally minimal.
\end{observation}

In fact, only the ``minimal" part of Observation~\ref{obs:sepsg:tm} is 
used, and the condition could be further weakened, as it was pointed out 
above. The reason for choosing hereditary  $h$-completeness as the 
condition, beyond its relatively simple definition, is its feature that 
every continuous homomorphic image of a hereditarily $h$-complete group 
is hereditarily $h$-complete. On the other hand, this property is not that 
far from $c$-compactness, because by \cite[2.16]{DikUsp}, {\slshape every 
hereditarily $h$-complete SIN group is $c$-compact}. ($G$ is {\em SIN} if 
its left and right uniformities coincide.) A second illustration for the 
phenomenon described in the Introduction is \cite[5.1]{DikUsp}, stating 
that {\slshape every locally compact $c$-compact  group is compact}. It is 
proved using Iwasawa's \cite{Iwas} Theorem:
{\slshape a connected locally compact group which has no closed subgroups 
(topologically) isomorphic to $\mathbb{R}$ is compact}. Dikranjan and 
Uspenskij's argument can be adjusted to hereditary $h$-completeness: since 
$\mathbb{R}$ is abelian and non-compact, it is not $h$-complete (see 
\cite[3.7]{DikUsp}), so no hereditarily $h$-complete group can contain it 
as a closed subgroup. Thus, one concludes:

\begin{observation} \label{obs:LC+hhc+conn:compact}
Every connected locally compact hereditarily $h$-complete group is 
compact.
\end{observation}

An obvious example is \cite[3.10]{DikUsp}, where instead of
$c$-compactness, the condition of ``all closed subgroups of $G$ are
$h$-complete" is imposed. As a corollary, one can observe that 
{\slshape every hereditarily $h$-complete soluble group is compact}, which 
is a slight generalization of \cite[3.12]{DikUsp}. We conclude this series 
of observations with one related to a result of the author. A group is 
called {\it maximally almost periodic} (or briefly, {\em MAP}) if it admits 
a continuous monomorphism $m: G \rightarrow K$ into a compact group $K$.
According to \cite[Cor.~7]{GL6}, {\slshape every $c$-compact MAP group is
compact}, but from \cite[Thm.~6]{GL6} one can also derive that {\slshape
every hereditarily $h$-complete MAP group is compact}.

The examples above show that so far no one seems to have exploited the
``extra" that $c$-compactness appears to have, compared to hereditary
$h$-completeness. This raises the following three problems:

\begin{problem} \label{prob:hhc:cc}
Is every hereditarily $h$-complete group $c$-compact?
\end{problem}

\begin{problem} \label{prob:wider}
Is there a class, wider than that of the SIN groups, in which the notions 
of hereditary $h$-completeness and $c$-compactness coincide?
\end{problem}

A universal algebraic approach might shed some more light on 
Problem~\ref{prob:hhc:cc}: For a class $\mathcal{W}$ of topological 
groups, put $\mathbf{P}(\mathcal{W})$ for the class of their 
(arbitrary) products, $\mathbf{H}(\mathcal{W})$ for the class of their 
continuous homomorphic images, and $\mathbf{\bar S}(\mathcal{W})$ for the 
class of the closed subgroups of groups from $\mathcal{W}$. A group $G$ is 
$h$-complete if $\mathbf{H}(G)$ consists of (Ra\u{\i}kov-)complete groups.
Thus, $G$ is hereditarily $h$-complete if each group in $\mathbf{\bar S}(G)$ 
is $h$-complete, in other words, if every group in 
$\mathbf{H}\mathbf{\bar S}(G)$ is complete. Since the product of any 
family of $h$-complete groups is $h$-complete (see \cite[2.13]{DikUsp}), 
if each group in $\mathbf{\bar S}(G)$ is $h$-complete, then so is every 
group in $\mathbf{P}\mathbf{\bar S}(G)$, and therefore 
every group in $\mathbf{H}\mathbf{P}\mathbf{\bar S}(G)$ is complete.
On the other hand, the product of any family of $c$-compact groups is 
$c$-compact (see \cite[2.8]{DikUsp} and \cite[4.4]{CleTho}), and 
$c$-compactness is a closed-hereditary property. Thus, if $G$ is 
$c$-compact, then every group in $\mathbf{\bar S}\mathbf{P}(G)$ is 
$c$-compact, and in particular $h$-complete; therefore, the groups in
$\mathbf{H}\mathbf{\bar S}\mathbf{P}(G)$ are complete. In most cases, one 
has $\mathbf{H}\mathbf{P}\mathbf{\bar S}(G) \subsetneq
\mathbf{H}\mathbf{\bar S}\mathbf{P}(G)$, which leads to a third problem:

\begin{problem} \label{prob:HSP}
Is hereditary $h$-completeness preserved under the formation of arbitrary 
products?
\end{problem}

Denote by I, II and III, respectively, the statements that the answers to
Problems I, II and III are affirmative. Implications \ref{prob:hhc:cc}
$\Longrightarrow$ \ref{prob:wider}, \ref{prob:HSP}, and \ref{prob:hhc:cc}
$\Longrightarrow$ \ref{prob:HSP} are trivial, but we are unable to say
anything further. It might be tempting to try to prove \ref{prob:HSP}
$\Longrightarrow$ \ref{prob:hhc:cc} (or even a stronger statement, that if
every group in $\mathbf{H}\mathbf{\bar S}\mathbf{P}(G)$ is complete, then
$G$ is $c$-compact), but we have a certain doubt about its truth. In any
case, it would be very interesting to characterize the structure of the
groups $G$ such that every group in $\mathbf{H}\mathbf{\bar
S}\mathbf{P}(G)$ is complete.

\bigskip

Solving the following problem is likely to be a key step towards 
understanding $c$-compactness and hereditary $h$-completeness:

\begin{problem} \label{prob:key}
Is there a hereditarily $h$-complete group that is not compact?
\end{problem}

If such a group $E$ exists, then either $E$ is $c$-compact, in which case 
it is a negative solution to the problem of $c$-compactness, or $E$ is not 
$c$-compact, and thus it provides a negative solution to
Problem~\ref{prob:hhc:cc}. On the other hand, if no such $E$ is found, 
then every hereditarily $h$-complete group is compact, and in particular
every $c$-compact group is compact. In the latter case, hereditary 
$h$-completeness, $c$-compactness and compactness coincide.

\section{Open map properties}
\label{sect:open}

A subset $F$ of a topological space $X$ is {\em $\omega$-closed} if 
$\bar C \subseteq F$ for every countable subset $C\subseteq F$.
A topological space $X$ has {\em countable tightness} if every 
$\omega$-closed subset is closed in $X$.

\begin{proposition} \label{prop:proper}
Let $(G,\mathcal T)$ be a hereditarily $h$-complete topological group, and 
let $\mathcal T^\prime$ be a coarser group topology on $G$. Then:

\begin{list}{{\rm (\alph{enumi})}}
{\usecounter{enumi}\setlength{\labelwidth}{25pt}\setlength{\topsep}{2pt}
\setlength{\itemsep}{0pt} \setlength{\leftmargin}{20pt}}

\item
every $\omega$-closed subset in $\mathcal T$ is $\omega$-closed in 
$\mathcal T^\prime$;

\item
$\mathcal T^\prime$ and $\mathcal T$ have the same compact subspaces.

\end{list}
\end{proposition}

\begin{proof}
Let $\iota: (G,\mathcal T) \rightarrow (G,\mathcal T^\prime)$ be the 
identity map. For a countable subset $D$, let 
$S=\overline{\langle D \rangle}$ be the closed separable subgroup 
generated by $D$. By Observation~\ref{obs:sepsg:tm}, the group $S$ is 
totally minimal, so $\iota\rest_S$ is a homeomorphism. In 
particular, $\iota(\bar D)=\overline{\iota(D)}$ and $\iota\rest_{\bar D}$ 
is a homeomorphism, because both $S$ and $\iota(S)$ are closed 
in the respective topologies.


(a) If $D$ is a countable subset of $\iota(F)$, then 
$\bar D = \iota(\overline{\iota^{-1}(D)})$. Thus, if $F$ is 
$\omega$-closed (in $\mathcal T$), then 
$\overline{\iota^{-1}(D)} \subseteq F$, and therefore 
$\bar D =\iota(\overline{\iota^{-1}(D)}) \subseteq \iota(F)$. 

(b) 
Let $K$ be a $\mathcal T^\prime$-compact subspace, and let $D$ be an 
countably infinite subset of $\iota^{-1}(K)$. Since $K$ is countably 
compact, it contains a limit point $y_0$ of $\iota(D)$. By the  foregoing 
discussion, $\iota\rest_{\bar D}$ is a homeomorphism, and thus 
$\iota^{-1}(y_0)$ is a limit point of $D$. Therefore, $\iota^{-1}(K)$ is 
countably compact. Clearly, $\iota^{-1}(K)$ is also closed; hence,
to complete the proof, we recall that a countably compact complete uniform 
space is compact \cite[p. 218, Ex.~21]{Wilan}.
\end{proof}

Before proving a more general result, we turn to what promises to be the
most simple case: discrete groups. Recall that a discrete group is
$c$-compact if and only if it is hereditarily $h$-complete
(\cite[5.3,2.16]{DikUsp}).

\begin{corollary}
Let $G$ be a discrete $c$-compact group. Then for every group topology
$\mathcal T$ on $G$:

\begin{list}{{\rm (\alph{enumi})}}
{\usecounter{enumi}\setlength{\labelwidth}{25pt}\setlength{\topsep}{2pt}
\setlength{\itemsep}{0pt} \setlength{\leftmargin}{20pt}}

\item
$\mathcal T$ is anticompact (i.e., every compact subset is finite);

\item
countably infinite subspaces are discrete in $\mathcal T$;

\item
every subset is $\omega$-closed in $\mathcal T$.
\qed
\end{list}
\end{corollary}

A map $f: X \rightarrow Y$ between Hausdorff spaces is {\em
$k$-continuous} if its restriction $f\rest_K$ to every compact
subspace $K$ of $X$ is continuous. Following Noble \cite{Noble}, a
topological group $G$ is called a {\em $k$-group} if every $k$-continuous
homomorphism $\varphi: G \rightarrow L$ into a topological group is
continuous. Every locally compact or sequential group is a $k$-group (they
are even {\em $k$-spaces}, see \cite[p. 152]{Engel6}). Noble showed in his 
papers that the class of $k$-groups is closed under the formation of 
quotients (in fact, continuous homomorphic images), arbitrary direct 
products, and open subgroups (see \cite[1.2,1.8]{Noble} and \cite[5.7]{Noble3}; 
the latter's statement is actually more general than what we quote here). 
This shows that the class of $k$-groups is a quite large one, and the 
only property it is missing to qualify for a variety is closedness 
under the formation of (arbitrary) subgroups. 

\pagebreak[3]

\begin{theorem} \label{thm:openmap}
Let $G$ be a hereditarily $h$-complete group, and let $f: G \rightarrow H$ 
be a continuous homomorphism onto a topological group $H$. 
If either
\begin{list}{{\rm (\alph{enumi})}}
{\usecounter{enumi}\setlength{\labelwidth}{25pt}\setlength{\topsep}{2pt}
\setlength{\itemsep}{0pt} \setlength{\leftmargin}{20pt}}

\item
$H$ has countable tightness, or

\item
$H$ is a $k$-group,
\end{list}

then $f$ is open.
\end{theorem}

\begin{proof}
Since hereditary $h$-completeness is preserved by quotients, by 
replacing $G$ with $G/\ker f$ and $f$ with the induced 
homomorphism, we may assume that $f$ is bijective. Set $g=f^{-1}$. 

(a) To show the continuity of $g$, let $F$ be a closed subset of $G$. Then 
$F$ is $\omega$-closed, and by Proposition~\ref{prop:proper}(a), 
$g^{-1}(F)=f(F)$ is $\omega$-closed in $H$. Since $H$ has countable 
tightness, this implies that $g^{-1}(F)$ is closed in $H$, and therefore 
$g$ is continuous.

(b)  Let $K \subseteq H$ be a compact subspace of $H$. By
Proposition~\ref{prop:proper}(b), the subset $g(K)=f^{-1}(K)$ is compact, 
so if $F$ is a closed subset of $g(K)$, then it is compact too. Thus,
$g^{-1}(F)=f(F)$ is compact, and in particular, $g^{-1}(F)$ is closed;
therefore, $g$ is $k$-continuous. Hence, $g$ is continuous, because $H$ is 
a $k$-group.
\end{proof}

The next corollary significantly generalizes \cite[3.2]{DikUsp}.

\begin{corollary} \label{cor:metr:open}
Let $G$ be a hereditarily $h$-complete topological group. Then every 
continuous homomorphism $f: G \rightarrow H$ onto a metrizable group $H$ 
is open. \qed
\end{corollary}

\begin{remark}
Corollary~\ref{cor:metr:open} might give the false impression of an Open 
Map Theorem that is free of Baire category arguments. However, this is not 
the case, because Proposition~\ref{prop:proper} relies on 
Observation~\ref{obs:sepsg:tm}, which in turn is a consequence of the 
classical Banach's Open Map Theorem for complete second countable 
topological groups.
\end{remark}

\section{Groups admitting a quasi-invariant basis}
\label{sect:quasi}

Following Kac \cite{Kac}, a group $G$ has a {\em quasi-invariant basis}
if for every neighborhood $U$ of the identity there exists a
countable family $\mathcal V$ of neighborhoods of the identity
such that for any $g\in G$ there exists $V \in \mathcal V$ such that
$g V g^{-1} \subset U$. In his paper, Kac proved the following 
two fundamental results on groups admitting a quasi-invariant basis.

\begin{fact}[{\cite{Kac}}]
\label{fact:Kac:embed}
A topological group can be embedded as a subgroup into a direct
product of metrizable groups if and only if it has a quasi-invariant
basis.
\end{fact}

\begin{fact}[{\cite{Kac}}] \label{fact:Kac:coarser}
A topological group $G$ with a quasi-invariant basis admits a coarser
metrizable group topology if and only if it has countable pseudocharacter.
\end{fact}

The result below is a structure theorem, which generalizes 
\cite[Thm.~5]{GL6} to the greatest possible extent, because every subgroup 
of a product of metrizable groups admits a quasi-invariant basis  
(Fact~\ref{fact:Kac:embed}).

\begin{theorem} \label{thm:QIB:struct}
Let $G$ be a hereditarily $h$-complete group  admitting a quasi-invariant 
basis. Then $G$ is the projective limit of its metrizable quotients.
\end{theorem}

\pagebreak[3]

\begin{proof}
Since $G$ has a quasi-invariant basis, by Fact~\ref{fact:Kac:embed}, it
embeds into a product $M=\prod\limits_{\alpha \in I} M_\alpha$ of metrizable 
groups. For $\pi_\alpha: G \rightarrow M_\alpha$ the restrictions of the 
canonical projections, we may assume that $\pi_\alpha$ is surjective, and 
thus by Corollary~\ref{cor:metr:open}, the $\pi_\alpha$ are open; in 
particular, each $M_\alpha$ is a (metrizable) quotient of $G$. Therefore, 
$G$ embeds into a product of its metrizable quotients, and since by 
adding  additional metrizable quotients to the product we cannot ruin the 
embedding property, we may assume that $M$ is the product of {\em all} the 
metrizable quotients of $G$.

If $G/N_1$ and $G/N_2$ are metrizable, then by Corollary~\ref{cor:metr:open}, 
the continuous monomorphism 
$G/N_1 \cap N_2 \rightarrow G/N_1 \times G/N_2$ is an embedding, 
so $G/N_1 \cap N_2$ is metrizable too. Thus, metrizable quotients of $G$ 
form a projective system. The image of $G$ is certainly contained in the 
projective limit of its metrizable quotients, and it 
is dense there. Since $G$ is $h$-complete, the statement follows.
\end{proof}

Dikranjan and Uspenskij \cite[3.3]{DikUsp} showed that {\slshape every
$h$-complete group with a countable network is totally minimal and
metrizable}. The statement concerning  metrizability can easily be
developed further:

\begin{proposition}
Every hereditarily $h$-complete group with a quasi-invariant basis that
has countable pseudocharacter is metrizable.
\end{proposition}

\begin{proof}
Let $(G,\mathcal T)$ be a group with the stated properties. By 
Fact~\ref{fact:Kac:coarser}, $G$ admits a coarser metrizable group 
topology $\mathcal T^\prime$. Applying Corollary~\ref{cor:metr:open} to 
the identity map $\iota: (G,\mathcal T) \rightarrow (G,\mathcal T^\prime)$,
one gets $\mathcal T = \mathcal T^\prime$, and thus $\mathcal T$ is 
metrizable.
\end{proof}

The difficulty with extending the first part of \cite[3.3]{DikUsp} 
(concerning minimality) is that given a hereditarily $h$-complete 
metrizable group $G$, the best we can say  (beyond  
Proposition~\ref{prop:proper}) is that each 
separable subgroup of $G$ is metrizable in every coarser topology. 

\begin{problem} \label{prob:met:tm}
Is every metrizable hereditarily $h$-complete group totally 
minimal?
\end{problem}

Since we do not know the answer to Problem~\ref{prob:met:tm}, we present a 
result relating it to the possibility of weakening the conditions of  
\cite[3.3]{DikUsp} in some sense. However, in order to do that, an 
auxiliary result is required, which turns out to be interesting on it own.
Its proof is modeled on the proof of \cite[3.4]{DikUsp}.

\begin{theorem} \label{thm:subdir:tm}
Let $G$ be a topological group and let $\mathcal{N}$ be a filter-base of 
$h$-complete normal subgroups of $G$. Suppose that $G$ naturally embeds 
into the product 
\[
P=\prod\limits_{N \in \mathcal{N}} G/N
\mbox{.}
\]
If each quotient $G/N$ is (totally) minimal, then $G$ is (totally)
minimal too.
\end{theorem}

\begin{remark}
The situation where $G$ is naturally embedded into a product of its 
quotients  is called a {\em subdirect representation}.
\end{remark}

\begin{proof}
Let $f: G \rightarrow H$ be a continuous surjective homomorphism, and let 
$U$ be an open neighborhood of the identity element in $G$. Since $G$ 
embeds into $P$, there exists a finite collection
$\{N_1, \ldots, N_k\} \subset \mathcal{N}$ such that 
$q^{-1}(V) \subseteq U$ for an open neighborhood $V$ of the identity  in 
the product $G/N_1 \times \cdots \cdots \times G/N_k$ (where $q$ is the 
diagonal of the respective projections). Since $\mathcal{N}$ is a 
filter base, there exists $N \in \mathcal{N}$ such that
$N \subseteq N_1 \cap \cdots \cap N_k$, and therefore without loss of 
generality we may assume that $U = p^{-1}(V)$ for an open neighborhood $V$ 
of the identity in $G/N$, where $p: G \rightarrow G/N$ is the canonical 
projection. 

Since $N$ is $h$-complete, $f(N)$ is closed in $H$, so $f$ induces a 
continuous surjective homomorphism $\bar f_N: G/N \rightarrow H/f(N)$. We 
have the following commutative diagram:
\[
\begin{CD}
G @>f>> H \\
@VpVV @VV\pi V \\
G/N @>\bar f_N>> H/f(N)
\end{CD}
\]
If $\bar f_N$ is  open, then so is the composite $\bar f_N \circ p$, and 
thus $\bar f (p(U)) = \pi(f(U))$ is open in $H/f(N)$. Therefore, 
$f(U)f(N)=f(UN) = f(U)$ is open in $H$ (one has $UN=U$, because 
$U = p^{-1}(V)$). Hence, $f$ is open if and only if $\bar f_N$ is open for 
every $N \in \mathcal{N}$.

One concludes that if each $G/N$ is totally minimal, then each $\bar f_N$,
being a continuous surjective homomorphism, is open, and therefore $G$ is
totally minimal. The argument for minimal groups is similar, because each
$\bar f_N$ is bijective when $f$ is so.
\end{proof}

Before returning to groups with quasi-invariant basis, we show how a
weaker version of \cite[7.3.9(b)]{DikProSto} follows from
Theorem~\ref{thm:subdir:tm}.

\begin{corollary} \label{cor:prod:tm}
Let $\{G_i\}_{i \in I}$ be a family of $h$-complete groups and put
$G = \prod\limits_{i \in I} G_i$. 

\begin{list}{{\rm (\alph{enumi})}}
{\usecounter{enumi}\setlength{\labelwidth}{25pt}\setlength{\topsep}{2pt}
\setlength{\itemsep}{0pt} \setlength{\leftmargin}{20pt}}

\item
If each $G_i$ is minimal, then $G$ is minimal too;

\item
If each $G_i$ is totally minimal, then $G$ is totally minimal too.
\end{list}
\end{corollary}

Part (a) of the Corollary is a new result, because 
\cite[7.3.9(b)]{DikProSto} deals only with products of totally minimal 
groups (and not minimal ones); on the other hand, part (b) is weaker, 
because in \cite[7.3.9(b)]{DikProSto} it suffices for the factors $G_i$ 
to have complete quotients, and they are not required to have complete 
homomorphic images. In order to prove Corollary~\ref{cor:prod:tm}, a weak 
version of \cite[(3)]{EDS} is needed:

\begin{fact} \label{fact:EDS}
The product of two (totally) minimal topological groups $G_1$ and $G_2$ is 
(totally) minimal if one of them is $h$-complete.
\end{fact}

\begin{proof}
For each finite subset $F \subseteq I$, put 
$N_F = \prod\limits_{i \in I\backslash F} G_i \times 
\prod\limits_{i \in  F} \{e\}$. By \cite[2.13]{DikUsp}, the product of any 
family of $h$-complete groups is $h$-complete, so each $N_F$ is 
$h$-complete. The quotient $G/N_F$ is topologically isomorphic to the 
finite product $\prod\limits_{i \in  F} G_i$, which is (totally) minimal 
by Fact~\ref{fact:EDS}, and the $N_F$ certainly form a filter-base.
Therefore, the conditions of Theorem~\ref{thm:subdir:tm} are fulfilled, 
and $G$ is (totally) minimal.
\end{proof}

A group $G$ is {\em perfectly totally minimal} if $G \times H$ is totally
minimal for every totally minimal group $H$ (\cite{Stoy}).

\begin{theorem} \label{thm:redu-QIB:tm}
The following statements are equivalent:

\begin{list}{{\rm (\roman{enumi})}}
{\usecounter{enumi}\setlength{\labelwidth}{25pt}\setlength{\topsep}{2pt}
\setlength{\itemsep}{0pt}}

\item
every metrizable hereditarily $h$-complete group is minimal;

\item
every hereditarily $h$-complete group with a quasi-invariant basis is
perfectly  totally minimal.
\end{list}
\end{theorem}

\begin{proof}
(i) $\Rightarrow$ (ii):
Since hereditarily $h$-completeness and metrizability 
are preserved by quotients, (i) implies that every metrizable 
hereditarily $h$-complete group is totally minimal.
By Theorem~\ref{thm:QIB:struct}, every hereditarily $h$-complete group $G$ 
with a quasi-invariant basis embeds into the product of its metrizable 
quotients, that are totally minimal by (i). Since $G$ is hereditarily 
$h$-complete, the kernel of each metrizable quotient is $h$-complete, and 
therefore, by Theorem~\ref{thm:subdir:tm}, $G$ is totally minimal. 
Hence, by Fact~\ref{fact:EDS}, $G$ is perfectly totally minimal, because 
it is $h$-complete.
\end{proof}

\section{Locally compact SIN groups}
\label{sect:SIN}

We recall that a group $G$ has {\em small invariant neighborhoods} (or
briefly, {\em $G$ is SIN}), if any neighborhood $U$ of $e \in G$ contains
an invariant neighborhood $V$ of $e$, i.e., a neighborhood $V$ such that
$g^{-1} V g = V$ for all $g \in G$. Equivalently, $G$ is SIN if its left
and right uniformities coincide. The interest in the class of SIN group in
the context of this paper arises from an aforesaid theorem of Dikranjan
and Uspenskij's \cite[2.16]{DikUsp}, stating that in this class,
$c$-compactness and hereditary $h$-completeness coincide.  It turns out
that locally compact SIN groups are closely related to discrete ones, and
as an illustration we start with two observations.

\begin{observation}
If every closed subgroup of a locally compact SIN group $G$ is totally 
minimal, then $G$ is $c$-compact.
\end{observation}

This statement was made concerning discrete groups in \cite[5.4]{DikUsp}.

\begin{proof}
Let $S$ be a closed subgroup of $G$, and let $f:S \rightarrow H$ be a 
continuous surjective homomorphism. Since $S$ is totally minimal, $f$ is 
open, and $H$ is a quotient of $S$. The group $S$, being a closed subgroup 
of $G$, is locally compact, and thus so is $H$; in particular, $H$ is 
complete. Therefore, $S$ is $h$-complete, and hence $G$ is hereditarily 
$h$-complete, which coincides with $c$-compactness, because $G$ is SIN.
\end{proof}

\begin{observation} \label{obs:sc+LC+SIN+cc:tm}
A $\sigma$-compact locally compact SIN group $G$ is $c$-compact if and 
only if every closed subgroup of $G$ is totally minimal.
\end{observation}

This mimics \cite[5.5]{DikUsp}, and it is an immediate consequence of the
previous observation and \cite[3.5]{DikUsp}. The condition of countability
in \cite[5.5]{DikUsp} was traded for $\sigma$-compactness, which is the
``right" concept for the context of LC groups.

\bigskip

Since locally compact connected $c$-compact groups are known to be compact
\cite[5.1]{DikUsp}, and compactness has the three-space property (i.e., if
$G/N$ and $N$ are compact for a normal subgroup $N$, then $G$ is compact
too), we turn first to totally disconnected groups.

\begin{proposition}
\label{prop:LC+TD+SIN:emb-discr}
Every locally compact totally disconnected SIN group is the projective 
limit of its discrete quotients with compact kernel.
\end{proposition}

\begin{remark}
In \cite[p. 68, 6.2]{GLPHD}, a different and more detailed proof of the 
Proposition is also available.
\end{remark}

\begin{proof}
In every locally compact totally disconnected group, the compact-open 
subgroups form a base at the identity. Since $G$ is also SIN, each 
compact-open subgroup contains a compact-open {\em normal} subgroup, so 
if we let ${\cal N}$ be the set of such subgroups, then $\mathcal{N}$ 
is a base at $e$. Thus, $G$ admits a continuous monomorphism $\nu$ into
$D=\prod\limits_{N \in {\cal N}} G/N$. The statement regarding the limit 
follows from a result established by  Weil in \cite[p. 25]{Weil}. 
\end{proof}

We are now ready to formulate a result of the same kind as
Theorem~\ref{thm:redu-QIB:tm}.

\begin{samepage}
\begin{theorem} 
The following statements are equivalent:

\begin{list}{{\rm (\roman{enumi})}}
{\usecounter{enumi}\setlength{\labelwidth}{25pt}\setlength{\topsep}{2pt}
\setlength{\itemsep}{0pt}}

\item
every discrete $c$-compact group is minimal;

\item
every locally compact $c$-compact group admitting small invariant 
neighborhoods is perfectly totally minimal.
\end{list}
\end{theorem}
\end{samepage}

\pagebreak[2]
In order to prove the Theorem, we need a corollary of Eberhardt, Dierolf
and Schwanengel:

\begin{fact}[{\cite[(7)]{EDS}}] \label{fact:3sp:tm}
If a topological group $G$ contains a compact normal subgroup $N$ such 
that $G/N$ is totally minimal, then $G$ is totally minimal.
\end{fact}

\begin{proof}
(i) $\Rightarrow$ (ii):
Let $G$ be a group described in (ii). By Fact~\ref{fact:EDS}, it suffices 
to show that $G$ is totally minimal, because it is $h$-complete. Set $N$ 
to be the connected component of the identity in $G$;  as a closed  
subgroup of $G$, $N$ is $c$-compact. Thus, $N$ is compact, because it is 
locally compact, connected and $c$-compact  
(Observation~\ref{obs:LC+hhc+conn:compact}). The quotient 
$G/N$ inherits all the aforesaid properties of $G$; furthermore, it is 
totally disconnected. By Fact~\ref{fact:3sp:tm}, it suffices to show that 
$G/N$ is totally minimal, so we may assume that $G$ is totally 
disconnected from the outset.

It follows from Proposition~\ref{prop:LC+TD+SIN:emb-discr} that $G$ embeds
into the product of its discrete quotients, where the kernels of the
quotients are the compact-open subgroups (certainly a filter-base), and in
particular they are $h$-complete. By (i), every discrete $c$-compact group
is minimal, and since quotients of discrete $c$-compact groups are again
discrete $c$-compact, (i) implies that they are actually totally minimal.
Thus, each discrete quotient of $G$ is totally minimal, and therefore, by
Theorem~\ref{thm:subdir:tm}, $G$ is totally minimal.
\end{proof}

Given a topological group $G$, the kernel of its Bohr-compactification 
$\kappa_G: G \rightarrow bG$ is called the {\em von Neumann radical} of 
$G$, and is denoted by $n(G)$. We say that $G$ is {\em minimally almost 
periodic} (or briefly, {\em m.a.p.}) if $n(G)=G$, or equivalently, if it 
has no non-trivial finite-dimensional unitary representations. Since 
$G$ is MAP if $n(G)$ is trivial, one can say that m.a.p. is the ``opposite" 
of MAP. We conclude the paper with a result similar in nature to 
\cite[Thm.~9]{GL6}.

\pagebreak[3]

\begin{theorem} \label{thm:red-SIN:cc}
The following statements are equivalent:

\begin{list}{{\rm (\roman{enumi})}}
{\usecounter{enumi}\setlength{\labelwidth}{30pt}\setlength{\topsep}{2pt}
\setlength{\itemsep}{0pt}}

\item
every countable discrete $c$-compact minimally almost periodic group is 
trivial;

\item
every countable discrete $c$-compact group is maximally almost periodic 
(and thus finite);

\item
every locally compact $c$-compact group admitting small invariant 
neighborhoods is compact.
\end{list}
\end{theorem}

\begin{proof}
(i) $\Rightarrow$ (ii):
The von Neumann radical of every $c$-compact group is m.a.p. 
(\cite[Cor.~8]{GL6}). Thus, by (i),  if $G$ is a countable discrete 
$c$-compact group, then $n(G)$ is trivial; in other words,  $G$ is MAP. 
The finiteness of $G$ follows from \cite[Cor.~7]{GL6}, but for the sake 
of completeness we show it explicitly here: $G$ is MAP, so
it admits a continuous monomorphism $m: G \rightarrow K$ into a compact 
group $K$. The map $m$ is an embedding, because $G$ is a discrete countable 
$c$-compact group, and as such it is totally minimal 
(Observation~\ref{obs:sc+LC+SIN+cc:tm}). Its image, $m(G)$, 
is closed in $K$, and thus compact, because $G$ is $h$-complete. 
Therefore,  $G$ is topologically isomorphic to the compact group $m(G)$.

(ii) $\Rightarrow$ (iii):
Let $G$ be a group as stated in (iii). The connected component $N$ of the 
identity in $G$ is compact (Observation~\ref{obs:LC+hhc+conn:compact}), so 
it suffices to show that $G/N$, which inherits all the aforesaid 
properties of $G$, is compact. Thus, we may assume that $G$ is totally 
disconnected from the outset.

If $S$ is a closed separable subgroup of $G$, then $S$ is $c$-compact, and 
its discrete quotients are countable. Since, by (ii), every discrete 
quotient of $S$ is finite, and $S$ is the projective limits of its 
discrete quotients (Proposition~\ref{prop:LC+TD+SIN:emb-discr}), $S$ is a 
{\em pro-finite} group; in particular, $S$ is compact. So every closed 
separable subgroup of $G$ is compact; thus, $G$ is precompact, and 
therefore, being complete, it is compact.
\end{proof}

It is not known whether (i) or (ii) is true. The only related result 
known is a negative one, due to Shelah \cite{Shelah}: under CH there 
exists an infinite discrete $h$-complete group. It is unknown whether such 
an example is available under ZFC. Even if such an example were available
under ZFC, it would not be, of course, a counterexample for (ii), because
it is neither countable nor $c$-compact.

Characterizing countable discrete MAP groups is a pure algebraic problem
that does not involve group topologies at all. On the other hand, a 
countable discrete group is $c$-compact if and only if its subgroups are 
totally minimal. By Theorem~\ref{thm:red-SIN:cc}, putting the two 
ingredients together can give a  solution to the problem of 
$c$-compactness in the locally compact SIN case.

An alternative approach could be to extensively study  countable 
discrete m.a.p. groups, which is again a purely algebraic task. Once a
characterization of such groups is obtained, it should no longer be 
difficult to check whether statement (i) in the Theorem holds.


\section*{Acknowledgments}

My deepest thanks go to my PhD supervisor, Professor Walter Tholen, for
introducing me to the problem of $c$-compactness. I would like to thank
him also for his attention to this paper.

I am grateful to Professor Dikran Dikranjan, my PhD external examiner, for
his incisive and helpful criticism, which has strongly influenced my
approach to topological algebra.

I am thankful to Professor Horst Herrlich for his attention to 
my work and his assistance in improving the presentation of the paper.

The constructive comments of the anonymous referee led to an improved
presentation of the results in this paper, for which I am grateful.

I wish to thank Professor Steve Watson, Professor Sakai Masami, and 
Professor Hans-Peter K\"unzi for the valuable discussions that were 
of great assistance in writing this paper.  I also appreciate the helpful 
comments of Professor Sasha Ravsky.

\bibliography{herhc}

\begin{samepage}

{\bigskip\bigskip\noindent 
FB 3 - Mathematik und Informatik\\
Universit\"{a}t Bremen\\ 
Bibliothekstrasse 1 \\
28359 Bremen \\
Germany

\nopagebreak
\bigskip\noindent{\em e-mail: lukacs@mathstat.yorku.ca} }
\end{samepage}

\end{document}